\newif\ifEMSstyle

\EMSstyletrue
\EMSstylefalse

\ifEMSstyle

\documentclass{book}
\usepackage[contrib,lang=british]{ems-book}

\else

\documentclass[12pt]{article}
\usepackage{amssymb,amsthm,mathrsfs,amsmath}
\usepackage[margin=2cm]{geometry}

\fi

\usepackage{mathrsfs}

\renewcommand{\geq}{\geqslant}
\renewcommand{\leq}{\leqslant}
\renewcommand{\ge}{\geqslant}
\renewcommand{\le}{\leqslant}

\newcommand\sym{\mathcal{S}}
\newcommand\PateC{\mathscr{C}}
\newcommand\norm[1]{\lVert#1\rVert}

\newcommand\eref[1]{$(\ref{#1})$}
\newcommand\sref[1]{\S$\ref{#1}$}

\newcommand\tref[1]{Theorem~$\ref{#1}$}
\newcommand\cjref[1]{Conjecture~$\ref{#1}$}
\newcommand\cyref[1]{Corollary~$\ref{#1}$}
\newcommand\fref[1]{Figure~$\ref{#1}$}

\newtheorem{thm}{Theorem}

\newtheorem{cor}{Corollary}
\newtheorem{conj}{Conjecture}

\newcommand\per{\operatorname{per}}

\newcommand\psd{{\cal H}}
\newcommand\corel{{\cal C}}
\newcommand\id{{\bf1}_n}
\newcommand\po{\preccurlyeq}

\newcommand{\pdc}{\cjref{cj:pdc}}
\newcommand{\pot}{\cjref{cj:pot}}

\begin{document}

\title{Lieb's permanental dominance conjecture}

\ifEMSstyle
\else

\author{Ian M. Wanless\\
  \small School of Mathematics\\[-0.75ex]
  \small Monash University, Australia\\
  \small \texttt{ian.wanless@monash.edu}
}

\date{}

\maketitle

\fi

\begin{abstract}
We survey the impact of Lieb's influential paper
``Proofs of some conjectures on permanents''
[{\em J.\ Math.\ Mech.\/} {\bf16} 1966, 127--134],
which introduced the famous permanental dominance conjecture.
This conjecture has defied all attacks for over half a century,
although a number of
related conjectures have recently been resolved.
\end{abstract}

\section{Introduction}

This is a survey article focusing on the legacy of Lieb's paper
\cite{Lie66} on permanents of matrices in $\psd_n$, the set of $n\times
n$ positive semi-definite Hermitian matrices.
Let $\sym_n$ denote the symmetric group on $\{1,2,\dots,n\}$ and
$\id$ denote the identity permutation in $\sym_n$. The {\em
permanent} of an $n\times n$ complex matrix $A=[a_{i,j}]$ is defined by
\[\per A=\sum_{\sigma\in \sym_n}\prod_{i=1}^n a_{i,\sigma(i)}.\]
More generally, 
if $G$ is a subgroup of $\sym_n$ and $\chi$ is any character of $G$
then the (normalised) {\em generalised matrix function} $f_\chi$ is defined by
\[f_\chi(A)=\frac{1}{\chi(\id)}\sum_{\sigma\in G}\chi(\sigma)\prod_{i=1}^na_{i,\sigma(i)}.\]
If $A\in\psd_n$ then $f_\chi(A)$ is a non-negative real number. If
$G=\sym_n$ and $\chi$ is irreducible then $f_\chi$ is called a (normalised)
{\em immanant}.  If $\chi$ is the principal/trivial character then
$f_\chi$ is the permanent, while if $\chi$ is the alternating
character then $f_\chi$ is the determinant. Taking $G$ to be the
trivial group yields the diagonal product
\[
h(A)=\prod_{i=1}^na_{i,i}.
\]
These three special examples of generalised matrix functions are related by
\begin{equation}\label{e:classical}
0\le\det A\le h(A)\le\per A
\end{equation}
for all $A\in\psd_n$. The second inequality was shown by Hadamard \cite{Had93}
and the last inequality is due to Marcus \cite{Mar63}, \cite{Mar64}.

Lieb's work \cite{Lie66} was motivated by results such as \eref{e:classical}.
He observed that ``the few
inequalities that are known for the permanent are suspiciously similar
to certain special cases of classical inequalities for the
determinants of such matrices -- the only difference being that the
direction of the inequality is reversed''.
Applying this principle to the classical result of Schur
\cite{Schur}, which states that $\det A\leq f_\chi(A)$ for
all $A\in\psd_n$, Lieb proposed:

\begin{conj}\label{cj:pdc} 
Let $G$ be a subgroup of $\sym_n$, and let $\chi$ be a character of $G$. Then
$\per A\geq f_\chi(A)$ for any $A\in\psd_n$.
\end{conj}

This conjecture subsequently became known as the ``permanental dominance
conjecture''. However, as noted by Zhang \cite{Zha16},
there is some confusion in the literature between this
name and the ``permanent on top conjecture''. We will use the latter name
for the related \pot\ below.


The special case of \cjref{cj:pdc} in which $\chi$ has degree 1 had
earlier been the subject of an open problem in \cite{MM65}, which was then
listed as Problem 2 in Minc's catalogue of open problems \cite{Minc1}.
\cjref{cj:pdc} was listed as Conjecture 42 in \cite{Minc3}.
As we shall see in \sref{s:relconj}, the permanental dominance conjecture
has motivated a lot of research, whilst gaining a degree of notoriety.
Despite considerable scrutiny, it remains open 56 years later.
The other key contribution from Lieb \cite{Lie66} was
the following result. An alternative proof of this theorem
was subsequently given by Djokovi\'c \cite{Djo69}.


\begin{thm}\label{t:Lieb}
  Let
  \begin{equation}\label{e:AinLiebthm}
  A=\left[
    \begin{array}{cc}
      B&C\\
      C^*&D\\
    \end{array}
    \right]\in\psd_n
  \end{equation}
  where $B$ and $D$ are $b\times b$ and $d\times d$ blocks, respectively.
  For any scalar $\lambda$, define 
  \[
  A_\lambda=\left[
    \begin{array}{cc}
      \lambda B&C\\
      C^*&D\\
    \end{array}
    \right].
  \]
  Now consider $P(\lambda)=\per A_\lambda$ as a polynomial in $\lambda$.
  All coefficients of this polynomial are real and nonnegative.
\end{thm}

Notice that $\per A=P(1)$ is the sum of the coefficients of $P(\lambda)$,
which therefore dominates any individual coefficient. In particular, it
dominates the coefficient of $\lambda^b$, which gives:

\begin{cor}\label{cy1}
  For $A$ in \eref{e:AinLiebthm} we have $\per A\ge (\per B)(\per D)$.
  Equality holds if and only if $A$ has a zero row or $C$ is the zero matrix.
\end{cor}

Also, if $b=d$ then $(\per C)(\per C^*)$ is the
constant term in $P(\lambda)$, which yields:

\begin{cor}\label{cy2}
  If $b=d$ in \eref{e:AinLiebthm} then
$
  \per A=P(1)\ge(\per B)(\per D)+(\per C)(\per C^*)$.
  Equality holds if and only if $A$ has a zero row or $C$ is the zero matrix.
\end{cor}

In his MathSciNet review of \cite{Lie66}, Marcus described the proofs
of these results as ``extremely ingenious and intricate arguments''.
Lieb also stated a determinantal analogue of \tref{t:Lieb}.

\bigskip

The following notation will be used throughout this paper.  The direct
sum of matrices $A$ and $B$ will be denoted $A\oplus B$, while their
tensor/Kronecker product will be denoted $A\otimes B$ and their Hadamard
(elementwise) product will be denoted $A\circ B$.
For any matrix $A$ the submatrix obtained by deleting row $i$ and
column $j$ from $A$ will be denoted $A(i|j)$. The Hermitian adjoint
(conjugate transpose) of $A$ will be denoted $A^*$, while
$\overline{A}$ will denote the conjugate of $A$.
We use $\corel_n$ to denote the subset of all matrices $[a_{i,j}]\in\psd_n$
that satisfy $a_{i,i}=1$ for all $1\le i\le n$. The matrices in $\corel_n$
are called \emph{correlation matrices}.

\section{Related conjectures}\label{s:relconj}

\newcommand\tick{_\checkmark}
\newcommand\wrong{_\times}
\newcommand\open{_\text{?}}

In this section we examine the relationships between a raft of
conjectures that are related to Lieb's permanental dominance
conjecture. To assist the reader to keep track of all these
conjectures, we summarise their relationships by the implications
pictured in \fref{f:relatedconj}. We show only the implications which
have historically been demonstrated in the literature prior to the
point at which a conjecture has been proved or refuted.  Of course,
once a conjecture has been resolved there are many new (uninteresting)
implications that could be added.  On the topic of resolution, the
current status of each conjecture is indicated in \fref{f:relatedconj}
by a subscript on the conjecture's number.  A $\tick$ indicates the
conjecture has been proved, a $\wrong$ shows that it has been
disproved and an $\open$ indicates that it remains open. These details
will gradually unfold in the historical account below.

\newcommand\cnj{\textrm{Conj }}
\begin{figure}[htb]
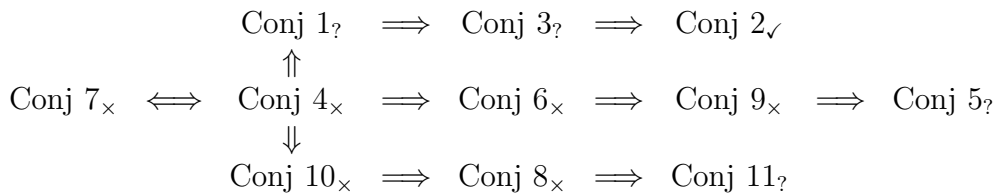

  \[\begin{matrix}
&    &\cnj\ref{cj:pdc}\open&\Longrightarrow&\cnj\ref{cj:Mar9}\open&\Longrightarrow&\cnj\ref{cj:marnew}\tick\\
&&\Uparrow&\\
\cnj\ref{cj:BS32}\wrong&\Longleftrightarrow&\cnj\ref{cj:pot}\wrong&\Longrightarrow&\cnj\ref{cj:BS38}\wrong&\Longrightarrow&\cnj\ref{cj:Beas}\wrong&\Longrightarrow&\cnj\ref{cj:Cho30}\open\\
&&\Downarrow&\\
&&\cnj\ref{cj:pate08}\wrong&\Longrightarrow&\cnj\ref{cj:BS40}\wrong&\Longrightarrow&\cnj\ref{cj:pate08b}\open\\
  \end{matrix}
  \]
\caption{\label{f:relatedconj}Relationship between various conjectures}
\end{figure}

We start by noting that Lieb's work was motivated by two
other conjectures.  Firstly, \cyref{cy1} solved the following
conjecture by Marcus and Newman, first published as Conjecture 8 in
\cite{MM65}. It is a permanental analogue of a classical result by
Fischer \cite{Fis07}, who showed that $\det A\le (\det B)(\det D)$
for $A$ partitioned as in \eref{e:AinLiebthm}.

\begin{conj}\label{cj:marnew}
  For $A$ in \eref{e:AinLiebthm} we have $\per A\ge (\per B)(\per D)$.
\end{conj}

Secondly, we have the following conjecture, which was
credited to Marcus as Conjecture 9 in \cite{MM65}.

\begin{conj}\label{cj:Mar9} 
Let $A\in\psd_{mk}$ be partitioned into
$k\times k$ blocks $A_{i,j}$, $i,j=1,2,\ldots,m$. Let $G$ be the
$m\times m$ matrix whose $(i,j)$ entry is
$\per(A_{i,j})$ for $i,j=1,2,\ldots,m$. Then
\begin{equation}\label{e:cj9}
\per A\ge \per G.
\end{equation}
If the $A_{i,i}$ are positive definite, then equality holds in
\eref{e:cj9} if and only if
\begin{equation*}
A=A_{11}\oplus A_{22}\oplus\cdots\oplus A_{mm}.
\end{equation*}
\end{conj}

It was noted in \cite{MM65} that \cjref{cj:Mar9} implies \cjref{cj:marnew}.
Also, Lieb~\cite{Lie66} showed that \cjref{cj:pdc} implies \cjref{cj:Mar9},
and his \cyref{cy2} proves the $m=2$ case. Pate \cite{Pate82} proved
that \cjref{cj:Mar9} holds when $A$ is real and $m$ is large enough
relative to $k$.

For any $A\in\psd_n$ we define the \emph{Schur power matrix} $\pi(A)$ to be
the $n!\times n!$ matrix whose $(\sigma,\tau)$ entry is
$\prod_{t=1}^na_{\sigma(t),\tau(t)}$, where $\sigma$ and $\tau$ run over all
permutations in $\sym_n$. 
The following conjecture was
introduced by Soules \cite{Soul66} and included as Conjecture 31 in
\cite{Minc1}.

\begin{conj}\label{cj:pot}  
  Let $A=[a_{i,j}]\in\psd_n$.
  Then $\per A$ is the maximum eigenvalue of $\pi(A)$.
\end{conj}

If $A\in\psd_n$ then
$\pi(A)$ is a principal submatrix of the $n$-fold
Kronecker product $\otimes^nA$, and hence $\pi(A)\in\psd_{n!}$.
It is easy to see that the row and column sums
of $\pi(A)$ equal $\per A$ and hence $\per A$ is an eigenvalue.
Schur \cite{Schur} (cf.~\cite{BS86}) implicitly showed
that $\det A$ is also an eigenvalue,
and indeed it is the lowest eigenvalue of $\pi(A)$.
So \pot\ is equivalent to the assertion that all eigenvalues of $\pi(A)$
lie in the interval $[\det A,\,\per A]$. It can be shown
(see Lemma 1 in \cite{Pate94b}) that $f_\chi(A)$ is
an eigenvalue of $\pi(A)$ for any character $\chi$ of a subgroup of $\sym_n$,
from which it follows that \pot\ implies \pdc.

Pate \cite{Pate89} proved a special case of \pot\ in order
to show that \pdc\ holds for immanants associated with two
part partitions. 
Soules \cite{Soul94} showed that if \pot\ is false for
real matrices then for the smallest order for which it fails there
must be a counterexample which is singular, has zero row sums and has
several other properties. Interestingly, as we will see shortly,
Drury \cite{Dru17} did eventually show that \pot\ fails for real matrices.

Showing similar intuition to that which motivated Lieb \cite{Lie66},
Chollet \cite{Chol} proposed the following permanental analogue of
Oppenheim's inequality for determinants:

\begin{conj}\label{cj:Cho30}
If $A,B\in\psd_n$ then
\begin{equation}\label{e:Chollet}
\per(A\circ B)\le(\per A)(\per B).
\end{equation}
\end{conj}

Chollet himself showed that it suffices to consider the case when
$B=\overline{A}$, when \eref{e:Chollet} reduces to
$\per(A\circ\overline{A})\le(\per A)^2$. By a standard scaling argument,
we also lose no generality by assuming that $A\in\corel_n$. The argument
goes like this. If $A=[a_{i,j}]$ has a zero row then both sides of
\eref{e:Chollet}
are zero and there is nothing to prove. So we may assume that $a_{i,i}>0$
for each $i$. Now define a diagonal matrix $D=[d_{i,j}]$ by
$d_{i,i}=a_{i,i}^{-1/2}$. Replace $A$ by $DAD$ and note that $DAD\in\corel_n$.
Also both $\per(A\circ\overline{A})$ and $(\per A)^2$ have been scaled by
the same factor, namely $h(A)^{-2}$, which justifies the claim.

Let $V=V(A)=\pi(A)\pi(A)^*=\pi(A)^2$. Note that for any $k$,
\begin{align*}
V_{k,k}&=\sum_{j=1}^{n!}\pi(A)_{k,j}\pi(A)^*_{j,k}
=\sum_{j=1}^{n!}|\pi(A)_{k,j}|^2=\sum_{\tau\in\sym_n}\Big|\prod_{i=1}^na_{i,\tau(i)}\Big|^2
=\sum_{\tau\in\sym_n}\prod_{i=1}^n|a_{i,\tau(i)}|^2\\
&=\per([|a_{i,j}|^2])=\per(A\circ\overline{A}).
\end{align*}
Also the sum of the elements in row $k$ of $V$ is
\[
\sum_{j=1}^{n!}V_{k,j}=\sum_{j=1}^{n!}\sum_{i=1}^{n!}\pi(A)_{k,i}\pi(A)_{i,j}
=\sum_{i=1}^{n!}\pi(A)_{k,i}\sum_{j=1}^{n!}\pi(A)_{i,j}
=\sum_{i=1}^{n!}\pi(A)_{k,i}\per A=(\per A)^2.
\]
So to prove \cjref{cj:Cho30} it suffices to show that $\pi(A)^2$ has row
sums that exceed its diagonal entries. Interestingly, if $A\in\corel_n$ then
$\pi(A)$ itself has the targetted property, since its row sums are 
$\per A$, which exceeds $h(A)=1=a_{k,k}$ for each $k$.

Zhang \cite{Zha13} initiated a new line of attack on \cjref{cj:Cho30}
by proving several properties of a ``maximising matrix'', namely
a matrix $X\in\corel_n$ that maximises $\per(A\circ X)$  for a given
$A\in\corel_n$.

Soon after Chollet's conjecture, Bapat and Sunder \cite{BS85}
proposed the following conjecture, which is stronger, given \eref{e:classical}.
Both conjectures are implied by \pot\ (see, for example, \cite{Merris87}).

\begin{conj}\label{cj:BS38}
If $A,B\in\psd_n$ then
$\per(A\circ B)\le (\per A)h(B)$.
\end{conj}

Using a similar scaling argument to the one we gave for
\cjref{cj:Cho30}, Zhang \cite{Zha89} noted that \cjref{cj:BS38} is
true if and only if it is true for all correlation matrices $A$ and
$B$. Zhang also showed that it is true if $A,B\in\corel_n$
and every off-diagonal entry of $B$ is equal to some fixed
$t$ in the interval $[0,1]$.

Bapat and Sunder \cite{BS86} proved \pot\ for $n\leq3$.  Marcus and
Sandy \cite{MS88} note that the $n=3$ case of \cjref{cj:Cho30} follows
immediately (see also Gregorac and Hentzel \cite{GH87}).
In the same paper, Bapat and Sunder \cite{BS86}
gave a reformulation of \pot\
which then appeared as Conjecture 32 in \cite{Minc3} as follows:

\begin{conj}\label{cj:BS32} 
Let $c$ be a complex valued function on $\sym_n$ satisfying
\[
\sum_{\sigma,\tau\in \sym_n}\overline{x(\tau)}\,c(\sigma\tau^{-1})\,x(\sigma)\ge0
\]
for all complex valued functions $x$ on $\sym_n$. Then
\[
c(\id)\per A\ge\sum_{\sigma\in \sym_n} c(\sigma)\prod_{i=1}^n a_{i,\sigma(i)}
\]
for any $A=[a_{i,j}]\in\psd_n$.
\end{conj}


Bapat and Sunder's paper \cite{BS86} also contained this conjecture:

\begin{conj}\label{cj:BS40} 
If $A$ is positive definite, then $\per A$ is the largest eigenvalue
of the matrix $\big[a_{i,j}\per A(i|j)\big]$.
\end{conj}

Much later, Beasley \cite{Beasley} made the following conjecture, which is
stronger than \cjref{cj:Cho30} and weaker than \cjref{cj:BS38}:

\begin{conj}\label{cj:Beas} 
If $A,B\in\psd_n$ then
\begin{equation}\label{e:cjBeas}
  \per(A\circ B)\le
  \max\big\{(\per A)h(B),\ (\per B)h(A)\big\}.
\end{equation}
\end{conj}

He proved this conjecture holds for $n\le3$ as well as noting that
\eref{e:cjBeas} is true if and only if it holds for all correlation
matrices (using a similar scaling argument to that seen above).

Let $A\in\psd_n$ and $1\le k\le n$. We define $\PateC_k(A)$ to be an $\binom
nk\times\binom nk$ matrix whose $\alpha,\beta$ entry is $\per
A[\alpha|\beta]\times \per A(\alpha|\beta)$. Here, the entries are
indexed by $\alpha,\beta\subset\{1,2,\dots,n\}$ with
$|\alpha|=|\beta|=k$ (we can impose lexicographic order on these row
and column indices, just for concreteness, but the order is not
important for our purposes). Also $A[\alpha|\beta]$ is the submatrix
induced by rows $\alpha$ and columns $\beta$, whereas $A(\alpha|\beta)$
is the submatrix formed by the removal of those rows and columns.
Pate \cite{Pate08} used \cyref{cy1} to
deduce that $\per A$ dominates the diagonal entries of $\PateC_k(A)$.  He
also showed that $\per A$ is an eigenvalue of $\PateC_k(A)$, and made the
following conjecture.
  
\begin{conj}\label{cj:pate08}
  Let $A\in\psd_n$ and $1\le k\le n$. 
  Then the largest eigenvalue of $\PateC_k(A)$ is $\per A$.
\end{conj}

The $k=1$ case of \cjref{cj:pate08} is precisely \cjref{cj:BS40}. Pate
showed that \cjref{cj:BS40} holds for a ``large subcone'' of $\psd_n$
and also for all nonnegative real matrices. Pate also raised the
challenge of proving the following consequence of \cjref{cj:BS40}:

\begin{conj}\label{cj:pate08b}
For $A\in\psd_n$,
\[a_{11}\per A(1|1)-a_{12}\per A(1|2)-a_{21}\per A(2|1)+a_{22}\per A(2|2)
\le2\per A.
\]
\end{conj}

\pot\ had achieved quite some prominence, in part because of the
number of consequences a proof would have (see \fref{f:relatedconj}).
However, in the last five years or so a cascade of counterexamples have
demolished most of those conjectures.  Shchesnovich \cite{Shc16} was
the first to disprove \pot, using a numerical search to find a rank 2
complex counterexample $A=W^*W$, where
\[
W=\left[\begin{array}{ccccc}
4-2i&2+3i&-4+4i&-3-4i&1\\
2+4i&3i&2+4i&3i&-5+7i\\
\end{array}\right].
\]
He found that $\per A=814\,016\,640$ is smaller than the largest eigenvalue
of $\pi(A)$, which is $320(2\,185\,775 +\sqrt{160\,600\,333\,345})$.

Drury \cite{Dru16} found a rank 2 complex $A\in\corel_7$
such that $A$ and $B=A^*$
provide a counterexample to \cjref{cj:BS38}. His $A=X^*X$ where
\[
X=\frac{1}{\sqrt2}\left[
  \begin{array}{ccccccc}
    \sqrt2&0&1&1&1&1&1\\
    0&\sqrt2&1&\omega&\omega^2&\omega^3&\omega^4\\
  \end{array}
  \right]
\]
and $\omega=e^{2\pi i/5}$ is a primitive fifth root of unity.
Note that $\per(A\circ B)=6185/128>45=\per A$.
Drury then followed that in \cite{Dru17} with a geometrically inspired
$16\times16$ real
counterexample to the same conjecture. He took a compound of
a regular dodecahedron
and a regular icosahedron inscribed together in the unit sphere.
He then chose one representative from each of the 16 pairs of antipodal vertices
and wrote its coordinates as one row in a $16\times 3$ matrix $Y$.
The real correlation matrix $A=YY^*$ provided a counterexample to
\cjref{cj:BS38} (with $B=A$). Drury also noted that $A$ gives a real
counterexample to \pot. He did not observe, but it is easy to see, that
both his counterexamples to \cjref{cj:BS38} are also counterexamples to
\cjref{cj:Beas}. Finally,
Drury \cite{Dru18} gave a rank 2 complex counterexample to \cjref{cj:BS40}. Let
$A=Z^*Z$, where \[
Z=\left[
\begin{array}{cccccccc}
-7+4i&9-3i&-6+2i&3+4i&7+6i&4-4i&i&5-8i\\
4-5i&1+4i&-8-2i&-7+4i&1-4i&1-8i&8-6i&1-3i\\
\end{array}
\right].
\]
Then $\per A=2\,977\,257\,622\,144\,118\,400$ and the
largest eigenvalue of $\big[a_{i,j}\per A(i|j)\big]$ is 
approximately $1.7\%$ larger at $\approx 3.028\times10^{18}$.

Tran \cite{Tra22} showed that \cjref{cj:pate08} is implied by \pot,
before developing the following (human checkable)
counterexample to both conjectures. Consider the rank 2 matrix
\[
H=\left[\begin{array}{ccccc}
3&1-2i&-1&1+2i&1\\
1+2i&3&1-2i&-1&1\\
-1&1+2i&3&1-2i&1\\
1-2i&-1&1+2i&3&1\\
1&1&1&1&1\\
\end{array}\right]\in\psd_5.
\]
Then $\per H=504$ and the spectrum of $\pi(H)$ is
$512^5$, $504^1$, $448^5$, $384^4$, $320^4$, $240^4$, $160^4$, $0^{93}$.
In particular,
$\per H$ is not the largest eigenvalue of $\pi(H)$, so $H$ is
a counterexample to \pot. Tran also showed that $H$
is a counterexample to \cjref{cj:pate08}.
Of course, on the basis of the implications in \fref{f:relatedconj},
\cite{Dru16}, \cite{Dru17} and \cite{Dru18} all provided
counterexamples to \pot. Interestingly, none of the counterexamples
mentioned so far appear to disprove \pdc, \cjref{cj:Cho30} or
\cjref{cj:pate08b}.

\bigskip

Both Drury \cite{Dru18} and Zhang \cite{Zha16} pose questions regarding
the smallest orders for which the disproved conjectures fail.
\pot\ is known to hold for $n\le3$ and fail for $n\ge5$.
The $n=4$ case remains open.
Shchesnovich \cite{Shc16} searched for counterexamples for order 4 but was unable to find any. Similarly, there is the question of the order of the smallest
real counterexample (currently the smallest known has order 16).
Another challenge is to search for new bounds on the eigenvalues of
$\pi(A)$ now that we know they can exceed $\per A$.

It is well-known that $\pi(A)$ is unitarily similar to a block
diagonal matrix where the blocks are indexed by the irreducible
representations of $\sym_n$. Hence each eigenvalue of $\pi(A)$ can be
associated with an irreducible character of $\sym_n$, and hence with a
partition of $n$. In the same way, any counterexamples
to \pot\ can also be associated with a particular partition of $n$. Drury
\cite{Dru18} calculated that the partitions associated with the
counterexamples in \cite{Shc16} and \cite{Dru18} are $(3,2)$ and $(7,1)$
respectively. He then conjectured that these are the only problematic
cases in the following sense:

\begin{conj}
  Suppose that the Ferrers diagram for a partition $\Pi$ of $n$
  does not contain the
  Ferrers diagram for $(3,2)$ nor $(7,1)$. Then for all $A\in\psd_n$,
  the largest eigenvalue of the block of $\pi(A)$ associated with
  $\Pi$ is at most $\per A$.
\end{conj}

To clarify what ``contains'' means in this last conjecture, the
Ferrers diagram for $\Pi$ contains $(3,2)$ if the first part of $\Pi$
is at least $3$ and the second part is at least $2$.
Drury \cite{Dru18} also conjectures that \cjref{cj:BS40} holds for
real matrices, and that \cjref{cj:BS38} holds for real matrices of
rank 2.

\section{Lieb's legacy}

In the previous section we examined many conjectures which had been
inspired either directly or indirectly by Lieb's paper
\cite{Lie66}. Those conjectures have undoubtedly prompted much
research of permanents of Hermitian matrices, even though many of the
conjectures are now known to be incorrect. In this section we examine
other aspects of the legacy of \cite{Lie66}.

We begin by briefly reviewing some of the highlights of progress
on \pdc. Any reader who is interested in further details is encouraged
to seek out the expositions by Merris \cite{Merris87},
\cite{Merris92}, James \cite{Jam87}, \cite{Jam92}, Pate
\cite{Pate94b}, \cite{Pate98}, \cite{Pate99},
Cheon and Wanless \cite{CW05}, Bapat \cite{Bap07}  and
Zhang \cite{Zha16}, among others.

Merris \cite{Merris83b} gave the following upper bound on generalised
matrix functions, in the spirit of \pdc, but known to be
weaker \cite{Merris87}. 

\begin{thm}
Let $G$ be a subgroup of $\sym_n$, and let $\chi$ be a character of $G$. Then
\[(h(A^n))^{1/n}\ge f_\chi(A)\]
for any $A\in\psd_n$.
\end{thm}

Most progress on the permanental dominance conjecture has been made on
its specialisation to immanants. For example, as a culmination of a
series of papers by Pate and others, we know \cite{Pate98} that
\pdc\ is true for all immanants when $n\le13$.  In comparison, the
general conjecture has only been shown for $n\le3$.  Tabata
\cite{Tab10}, \cite{Tab15} examined the $n=3$ case in detail and
showed strict inequality holds when $\chi$ is not the principal character.
Interestingly, James \cite{Jam92} discovered that the following
matrix in $\psd_4$
\[
\left[\begin{matrix}
\sqrt{3}&i&i&-{i}\\
-{i}&\sqrt{3}&i&i\\
-{i}&-{i}&\sqrt{3}&-{i}\\
i&-{i}&i&\sqrt{3}\\
  \end{matrix}\right]
\]
achieves equality in \cjref{cj:pdc}
when $G$ is the alternating group $A_4$ and $\chi$ is the
character of that group that satisfies $\chi((12)(34))=1$ and
$\chi((234))=e^{2\pi i/3}$. It follows that for any $n\ge4$ there
is at least one non-principal character which achieves equality
in \cjref{cj:pdc} for some $A\in\psd_n$.

To describe the results which have been proved on immanants we define
a partial order $\po$ on the set of partitions of an integer $n$.  Let
$\lambda$ and $\mu$ be two such partitions and let $\chi$ and $\chi'$
be the characters associated with $\lambda$ and $\mu$ respectively by
the well known bijection between partitions of $n$ and irreducible
characters of $\sym_n$.  By $\lambda\po\mu$ we will mean
that $f_\chi(H)\leq f_{\chi'}(H)$ for all $H\in\psd_n$.
Let $J_n\in\psd_n$ be the matrix in which ever entry is $1$. 
Note that, $f_{(3,1)}(J_2\oplus J_2)<f_{(2,2)}(J_2\oplus J_2)$
so it is not true that $(2,2)\po(3,1)$. Similarly,
$f_{(3,1)}(J_3\oplus J_1)>f_{(2,2)}(J_3\oplus J_1)$
so that $(3,1)\po(2,2)$ also fails.
Hence $\po$ is not a total order. 
James \cite{Jam92} notes that using similar
tests involving direct products of blocks $J_i$ for various values
of $i$, it can be shown that
for two partitions $\lambda,\mu$ of $n$ to satisfy $\lambda\po\mu$
it is necessary but not sufficient that $\mu$ majorises~$\lambda$.


The result of Schur \cite{Schur} implies that $(1^n)\po\lambda$ for
all partitions $\lambda$ of $n$.  The specialisation of the
permanental dominance conjecture to immanants asserts that for
all~$\lambda$,
\begin{equation}\label{e:pdc}
\lambda\po(n).
\end{equation}

The following beautiful theorem of Heyfron \cite{Hey88} shows that the
``single-hook immanants'' are neatly ordered between the determinant
and permanent.

\begin{thm}\label{t:Heyfron}
\begin{equation}\label{e:hook}
(1^n)\po(2,1^{n-2})\po(3,1^{n-3})\po\cdots\po(n-1,1)\po(n).
\end{equation}
\end{thm}

Heyfron's theorem confirmed a conjecture originally made by Merris
\cite{Merris83} (note that Merris \cite{Merris87} attributes the
conjecture to himself and Pierce).  A number of partial results
in this direction had been obtained by Merris and Watkins \cite{MW85} and
Johnson and Pierce \cite{JP87}, \cite{JP88}, for example.

Stembridge \cite{Stem91}, \cite{Stem92} notes that the analogue of the
permanental dominance conjecture is trivially true for totally
positive matrices (real matrices with non-negative minors). He
considers to what extent analogues of \tref{t:Heyfron} hold for this
class of matrices. He also considers inequalities for immanants of
so called Jacobi-Trudi matrices, which are closely connected with the
combinatorics of symmetric functions.  Haiman \cite{Hai93} develops
these ideas and makes connections with Kazhdan-Lusztig theory
and characters of Hecke algebras.

Chan and Lam \cite{CL97}, \cite{CL98} sharpened the inequalities in
\eref{e:hook} in the case of matrices which are the Laplacians
of trees. More recently, Nagar and Sivasubramanian \cite{NS17} 
generalised Chan and Lam's work to q-Laplacians.

A general scheme for obtaining inequalities
involving immanants is described by Pate in \cite{Pate97} and
many such inequalities can be found throughout his papers.
For example, in \cite{Pate99} he showed the following results for positive
integers $n$, $p$ and $k$.  If $k\geq2$ and $n\geq p+k-2$ then
$(n+p-i,n^k,i)\po(n+p,n^k)$ for $1\leq i\leq p$.  On the other hand if
$p\geq n+k-1$ then $(n+p-i,n^k,i)\po(n+p,n^k)$ whenever $p/2\leq i\leq
n$.  In the same paper he obtained the following asymptotic result.
For positive integers $k$ and $s$ there exists an integer $N_{k,s}$
such that for all $n\geq N_{k,s}$,
\[(n+s,n^k)\po(2n+s,n^{k-1})\po(3n+s,n^{k-2})\po\cdots\po(kn+n+s).
\]

\bigskip

\cyref{cy1} has proven useful when deriving special cases of the
permanental dominance conjecture (see e.g. Merris and Watkins
\cite{MW85}).  For example, it immediately proves that \pdc\ holds when
$G$ is a Young subgroup of $\sym_n$ (that is, a direct product of symmetric
groups) and $\chi$ is the trivial character on $G$.

Next we turn our attention to other applications of \tref{t:Lieb} and its
corollaries.
Strengthening \eref{e:classical},
Marcus and Soules \cite{MS67} found upper and lower bounds on $\per A-h(A)$
and $h(a)-\det A$ for $A\in\psd_n$. They also
proved this stronger form of \cyref{cy1}:

\begin{thm}
For $A$ in \eref{e:AinLiebthm} we have $\per A\ge (\per B)(\per
D)+\mu^{n-2}\norm{C}^2$, where $\mu$ is the smallest eigenvalue of $A$
and $\norm{\cdot}$ denotes the Euclidean norm.
\end{thm}

Grone and Pierce \cite{GP88} used \cyref{cy1} to prove the following
bound for permanents of correlation matrices.

\begin{thm}\label{t:GP}
  Let $A=[a_{i,j}]\in\corel_n$. Then $\per A\ge\frac1n\norm{A}$,
  with equality if and only if $n=2$ or $A$ is an identity matrix or $A$ is unitarily similar to
  \[
  Y_3=
  \left[
    \begin{array}{ccc}
      1&-\frac12&-\frac12\\[1ex]
      -\frac12&1&-\frac12\\[1ex]
      -\frac12&-\frac12&1\\
    \end{array}
    \right].
  \]
\end{thm}

From \tref{t:GP} they also deduced the following three corollaries, each of
which had been conjectured by Grone and Merris \cite{GM87}.

\begin{cor}\label{cy:hAA}
  If $A\in\corel_n$ then $\per A\ge(h(A^2))^{1/n}$.
\end{cor}

\begin{cor}\label{cy:mixdetper}
  If $A\in\psd_n$ then $(n-1)\per A + \det A\ge n\,h(A)$.
\end{cor}

\begin{cor}\label{cy:singcorrel}
  If $A\in\corel_n$ is singular then $\per A\ge n/(n-1)$.
\end{cor}

The following result of Pate~\cite{Pate03} was motivated by
\cyref{cy:singcorrel} (cf.~\eref{e:classical}).

\begin{thm}
  If $A\in\psd_n$ is real, and has rank at most $2$, then
  \[\per A\ge\frac{n!}{2^{n-1}}h(A).\]
\end{thm}

Frenkel \cite{Fre08} derives analogues of \tref{t:Lieb} for Pfaffians
and Hafnians, then applies his results to something called the
real linear polarisation constant problem.
He also remarks that in order to derive \tref{t:Lieb}
it is not necessary to assume that $A$ is positive semi-definite;
it suffices for both the diagonal blocks $B$ and $D$ to be
positive semi-definite. 

In a follow-up paper,
Frenkel \cite{Fre10} obtained an analogue of \cyref{cy1} for
$\alpha$-permanents. The $\alpha$-permanent is defined by
\begin{equation}\label{e:alphper}
\per{}_{\alpha}(A)=\sum_{\sigma\in\sym_n}\alpha^{\nu(\sigma)}\prod_{i=1}^na_{i,\sigma(i)},
\end{equation}
where $\nu(\sigma)$ is the number of disjoint
cycles of the permutation $\sigma$. Note that the 1-permanent is just the
permanent. Also, the $(-1)$-permanent is equal to the determinant, up to a
factor of $(-1)^n$. Frenkel's Theorem is this:

\begin{thm}
Suppose that $\alpha$ is either a nonnegative integer or a real number with
$\alpha\ge n-1$. Then, for $A\in\psd_n$ as in 
\eref{e:AinLiebthm} we have $\per_\alpha A\ge (\per_\alpha B)(\per_\alpha D)$
and $0\le(-1)^n\per_{-\alpha} A\le (-1)^n(\per_{-\alpha} B)(\per_{-\alpha} D)$.
\end{thm}

Another application of \cyref{cy1} to $\alpha$-permanents
is given by Eisenbaum \cite{Eis12} in her study of the stochastic
comparisons between $\alpha$-permanental point processes.
Shirai~\cite{Shi07} studies a function closely related to
\eref{e:alphper} and considers for what values of $\alpha$ it
is positive on $\psd_n$.

\bigskip

For a short paper, \cite{Lie66} has proved an inspiration for 
investigations in many different directions.

 
  \let\oldthebibliography=\thebibliography
  \let\endoldthebibliography=\endthebibliography
  \renewenvironment{thebibliography}[1]{%
    \begin{oldthebibliography}{#1}%
      \setlength{\parskip}{0.2ex}%
      \setlength{\itemsep}{0.2ex}%
  }%
  {%
    \end{oldthebibliography}%
  }

\end{document}